\documentclass[graybox,leqno]{svmult}

\usepackage{makeidx}         
\usepackage{graphicx}        
\usepackage{multicol}        
\usepackage[bottom]{footmisc}

\usepackage[english]{babel}
\usepackage{amssymb,amstext,amsmath}

\usepackage{color}
\definecolor{darkgreen}{cmyk}{1,0,1,.2}
\definecolor{m}{rgb}{1,0.1,1}

\usepackage{soul}

\makeindex             

\begin{document}

\title*{Some examples of first exit times}
\author{Jes\'us Antonio \'Alvarez L\'opez\inst{1}\and
Alberto Candel\inst{2}}
\institute{\inst{1}Departamento de Xeometr\'ia e Topolox\'ia,
  Facultade de Matem\'aticas, Universidade de Santiago de Compostela,
  15706 Santiago de Compostela, Spain; \texttt{jesus.alvarez@usc.es}
   \\
   \inst{2}Department of Mathematics, California State University,
  Northridge, CA 91330, U.S.A.; \texttt{alberto.candel@csun.edu}}
%
%
\maketitle

\abstract{The purpose of this article is to compute the expected first
  exit times of Brownian motion from a variety of domains in the
  Euclidean plane and in the hyperbolic plane.}

\section{Introduction}

The theory of Brownian motion on Riemannian manifolds allows
probabilistic interpretations of solutions to second order
differential equations on them via the so called Dynkin formula. One
example of such equation and solution is the following: if $D$ is a
regular domain on a Riemannian manifold $M$ with attending Laplace
operator $\triangle$, then the expected value of the ``first exit time''
from \(D\) for Brownian paths in the Wiener space of \(M\), if finite, is
the minimal solution to the differential equation
\[\triangle f \equiv -1\]
on $D$, with \(f>0\) on \(D\) and $f\equiv 0$ on $\partial
D$. The number \(\rho(D)=4\int_D f\) is called the torsional rigidity
of \(D\), a sort of isoperimetric constant whose study originated with
Saint-Venant memoir~\cite{SaintVenant}.

Ghys~\cite{Ghys95} gave a spectacular application of first exit times to the
topology and dynamics of foliated spaces. The examples and calculations
presented here are the base of some exercises in \cite[Exercises
2.7.14, C.9.4 and C.9.5]{CandelConlon04}, and were motivated by a
discussion on Ghys theorem.

\section{Generalities on first exit times and on harmonic functions}\label{section:general}

More details about first exit times and probabilistic solutions to differential
equations can be found in Dynkin~\cite{Dynkin65}; the theory needed
for what follows is detailed in~\cite[Appendix C]{CandelConlon04}.

Consider a Riemannian manifold, \(M\), with attending Laplacian
\(\triangle\). These data permit to construct Brownian motion on
\(M\), a continuous-time stochastic process taking place in the space
of continuous paths \(\omega:[0,\infty)\to M\) that is regulated by a set
of probability measures \(\{P_x\mid x\in M\}\) (with \(P_x\) supported
on paths \(\{ \omega(0)=x\}\)) which are constructed via the heat
kernel density of the Laplacian \(\triangle\).

Solutions to a variety of differential equations on \(M\) admit
probabilistic interpretations via Brownian motion. One such example is
the following. Let \(D\subset M\) be a regular domain (a connected
open set with piecewise smooth boundary), and consider the first exit time
\(T_D\) from \(D\), the function on paths given by
\(T_D(\omega)=\inf\{ t>0\mid \omega(t)\notin D\}\) (with the standard
convention that the infimum of the empty set is \(\infty\)).

The expected first exit with respect to the Brownian 
measures \(\{P_x\}\) defines a function \(E_\bullet[T_D]: x\mapsto
E_x[T_D]=\int T_D(\omega) \cdot P_x(\omega)\) which is \(0\) for all
\(x\notin D\), and which is either \(E_x[T_D]<\infty\) for all \(x\in
D\), or \(E_x[T_D]\equiv \infty\) for all \(x\in D\).

If \(D\) is a relatively compact domain, then \(E_x[T_D]<\infty\) for
all \(x\in D\). In this case, Dynkin's formula shows that this function
is a solution to the differential equation problem (Saint-Venant
problem)
\begin{svgraybox}
\begin{equation}\label{equation1} 
\left \{
\begin{array}{rll}
\triangle f &\equiv -1 &\quad\text{on \( D\)}, \\
 f&>0 &\quad\text{on \(D\), and} \\ 
f&\equiv 0 &\quad\text{on \(\partial D\)}.
\end{array}
\right.
\end{equation}
\end{svgraybox}

Because of the maximum modulus principle for harmonic functions (to
the effect that a function that is harmonic on a relatively compact
domain and continuous on its closure must attain its extreme values on
the boundary of the domain), the solution to the differential
equation~(\ref{equation1}) is unique: the difference of two solutions
is harmonic and equal to \(0\) on \(\partial D\), so it must be \(0\)
on all of \(D\). Therefore, if \(D\) is relatively compact, the
expected first exit time from \(D\) is the unique solution to the
differential equation~(\ref{equation1}).
  
If \(D\) is not relatively compact, then \(E_\bullet[T_D]\) may or may not
be a finite function. At any rate, there is an increasing
sequence of relatively compact domains \(D_1\subset D_2 \subset \ldots
\subset D\) that exhaust \(D\). The first exit time functions, \(T_{n}\), from
\(D_n\) increase pointwise to the first exit time function \(T_D\) and thus
the monotone convergence theorem implies that 
the sequence of expected first exit times \(E_\bullet[T_n]\) 
increases  to the expected first exit time \(E_\bullet[T_D]\).  
 
Furthermore, if \(E_x[T_D]<\infty\) for one \(x\in D\), then
\(E_\bullet[T_D]<\infty\) everywhere and is a solution to
Equation~(\ref{equation1}). In fact, if the expected first exit time
\(f=  E_\bullet[T_D]<\infty\), then \(f\) is the minimal solution to
that equation on \(D\).  Indeed, from the above paragraph, you infer
that \(f=E_\bullet[T_D]\) is given by
\[f=\sup_{B\subset D} f_B\] where \(\{ B\subset D\} \) is the set of
relatively compact regular domains contained in \(D\), and where
\(f_B=E_\bullet[T_B]\) is the expected first exit time from \(B\). If
\(g\) is any positive function on \(D\) such that \(\triangle g=-1\)
and \(g\equiv 0\) on \(\partial D\), then \(g\ge f\) because, if that was
not the case, then \(g<f\) on an open subset of \(D\), and so it
follows from the definition \(f=\sup_B f_B\) that there exists a
relatively compact domain \(B\subset D\) where \(g< f_B\). Then
\(f_B-g\) is harmonic and \(>0\) on \(B\) but \(\le 0\) on
\(\partial B\), in contradiction to the maximum modulus principle.
 
Moving on to a brief review of harmonic functions, besides the already
mentioned maximum modulus principle, two other well-known facts will
be repeatedly used below.  Both concern harmonic functions on domains
in \(\mathbf{R}^2\) endowed with Riemannian metrics conformal to the
standard metric, that is, of the form \(\varphi(dx\otimes dx+dy\otimes
dy)\). The second fact is then that, since the Laplacian for this
metric is \(\triangle_\varphi
u=(1/\varphi)\left(u_{xx}+u_{yy}\right)\), a function \(u\) is
harmonic on a domain of this type if and only if \(u\) is harmonic in
the classical sense that \(u_{xx}+u_{yy}=0\) all throughout the
domain. Because of this, if \(\phi\) is a holomorphic function with
range in the domain of the harmonic function \(u\), then the composite
\(u\circ \phi\) is harmonic in the domain of \(\phi\), as is easily
verified via the chain rule, utilizing the harmonicity of \(u\) and
the Cauchy-Riemann equations for \(\phi\).

The third fact about harmonic functions is deeper and concerns their
integral representations. To each function \(u\ge 0\) that is harmonic
on the right half plane \(\{x>0\}\) there corresponds a measure
\(\mu\) on the line \(\Im=\{(0,t)\mid -\infty<t<\infty\}\) and a
constant \(C\ge 0\) so that
\[ u(x,y)=C x + \int_\Im \frac{x}{x^2+(y-t)^2} \cdot \mu(t),\]
for all \((x,y)\) with \(x>0\). 

In particular, if \(u\ge0\) is harmonic and extends continuously by
\(0\) to all but finitely many points \(t_1, t_2, \ldots, t_n\) on the
line \(x=0\), then the measure \(\mu\) is supported on the set
\(\{t_k\}\), and \(u\) may be expressed as
\begin{equation}\label{PoissonKernelRep}
 u(x,y)=  C_\infty x + \sum_{k=1}^n \frac{C_k x}{x^2+(y-t_k)^2}, 
\end{equation}
for some constants \(C_\infty, C_1, C_2, \ldots, C_n\ge 0\).

\section{Domains in the Euclidean plane}
\label{sec:euclidean}

In Cartesian coordinates $(x,y)\in \mathbf{R}^2$, the Euclidean metric is
\( dx\otimes dx + dy\otimes dy\), 
and the Laplacian is given by
\[ \triangle_e f= f_{xx} +f_{yy}. \] 
In  polar coordinates $(r,\theta)$, the Laplacian is given by 
\begin{equation}\label{polar laplacian}
\triangle_e f = f_{rr} +\frac{1}{r} f_r +\frac{1}{r^2} f_{\theta
\theta}.\end{equation}

\subsection{Domain bounded by an ellipse}
\label{subsec:ellipse}

Let \(D\) be the domain enclosed by an ellipse in the Euclidean plane
with axes of lengths \(a, b>0\) and center \((h,k)\). Up to isometry
(a rotation), \(D\) consists of all points \((x,y)\) such that
\(\dfrac{(x-h)^2}{a^2}+\dfrac{(y-k)^2}{b^2}<1\).

\begin{svgraybox}
The expected first exit time from \(D\) is given by the function
\begin{equation}\label{exit:ellipse}
f(x,y)=\dfrac{a^2b^2}{2a^2+2b^2}
\left(1-\frac{(x-h)^2}{a^2}-\frac{(y-k)^2}{b^2}\right).
\end{equation}
\end{svgraybox}

Indeed, \(f\) is positive on \(D\), identically \(0\) on the ellipse
\(\partial D\), and satisfies the differential equation \(\triangle
f=-1\) on \(D\). Therefore \(E_{(x,y)}[T_D]=f(x,y)\) because, \(D\) being
relatively compact, Equation~(\ref{equation1}) has exactly one
solution.

\subsection{Domain bounded by a parabola}\label{section:parabola}

Up to isometry, a parabola has an equation of the form \(y^2={4p x}\)
(focus at \((p,0)\) and focal distance \(p\)), and a convex domain
bounded by a parabola is isometric to the domain, \(D\), consisting of
all \((x,y)\in \mathbf{R}^2\) such that \(4px> y^2\).
\begin{svgraybox} The expected first exit time from \(D\) is
\begin{equation}\label{exit:parabola}
E_{(x,y)}[T_D]= 2px-\dfrac{y^2}{2},\end{equation}
for all \((x,y)\in D\).
\end{svgraybox}

It is plain that \(f(x,y)=2px-y^2/2\) is a solution to
equation~(\ref{equation1}) on \(D\), but to prove that the expected
first exit time \(E_{(x,y)}[T_D]=f(x,y)\) on \(D\) requires some extra work
because \(D\) is not relatively compact. Let \(D_n\), \(n >[p]\),
denote the domain enclosed by the ellipse with foci at \((p,0)\) and
\((2n-p,0)\) and eccentricity \(e=1-p/n\).  An equation for  this
ellipse is
\[ \frac{(x-n)^2}{n^2} + \frac{y^2}{2p(n-p)}=1.\] For \(n>[p]\),
the domains \(D_n\subset D_{n+1}\subset\ldots \subset D\) increase to
the domain \(D\), and so  the first exit times functions \(T_n\) of
\(D_n\) increase pointwise to the first exit time function \(T_D\). By the
dominated convergence theorem, \(E_{(x,y)}[T_n]\) converges to
\(E_{(x,y)}[T_D]\).  The expected first exit time
from \(D_n\) was shown to be (\ref{exit:ellipse})
\begin{align*} E_{(x,y)}[T_n] & = \frac{n^2p(n-p)}{2p(n-p)+n^2}
  \left( 1-  \frac{(x-n)^2}{n^2} - \frac{y^2}{2p(n-p)}\right)\\
& =-\frac{p(n-p)}{2p(n-p)+n^2} x^2 +
 \frac{2np(n-p)}{2p(n-p)+n^2} x
  - \frac{n^2}{4p(n-p)+2n^2} y^2.
\end{align*}
It follows immediately that \( E_{(x,y)}[T_n] \to 2px- y^2/2\),
as \(n\to \infty\), which is the
expression for the solution to equation~(\ref{equation1}) shown
at~(\ref{exit:parabola}),  uniformly on compact subsets of \(D\).

Any solution to Equation~(\ref{equation1}) is of the form
\(E_\bullet[T_D] +u\), where \(u\) is harmonic and \(\ge0\) on \(D\)
and identically \(0\) on \(\partial D\).  The function \(\phi:z\mapsto \cosh
\dfrac{\pi}{2}\sqrt{\dfrac{z}{p}-1}\) is a conformal representation of
\(D\) onto the right half plane \(\Re z>0\). If \(u\) is harmonic and
\(\ge0\) on \(D\) and identically \(0\) on \(\partial D\), then
\(u\circ \phi^{-1}\) is harmonic and \(\ge0\) on \(\Re z>0\) and
identically \(0\) on \(\Re z=0\), so, by~(\ref{PoissonKernelRep}),
\(u\circ\phi^{-1}(w)= C\Re w\), for some \(C\ge 0\),
or \(u(z)=C \phi(z)\) after the  switch \(w=\phi(z)\). 

\begin{svgraybox} That is,
any solution to Equation~(\ref{equation1}) on \(D\) is given by (using
complex coordinates \(z=x+y i\)):
\begin{equation}\label{equation1:parabola}
z\mapsto 2 p \Re(z) -\frac{\Im z^2}{2}+ C \Re \left(\cosh
  \dfrac{\pi}{2}\sqrt{\dfrac{z}{p}-1}\right),
\end{equation}
for some constant \(C\ge0\). 
\end{svgraybox}

 \subsection{Domain between two concentric circles}

 Let \(D\) be a domain bounded by two concentric circles of radii
 \(a<b\). Rotations about the common center of the circles are
 isometries that leave \(D\) invariant. Therefore the expected first
 exit time, \(f=E_\bullet[T_D]\), from \(D\), in polar coordinates
 \((r,\theta)\) about its center is a function of \(r\) only, and so,
 in those coordinates, the equation \(\triangle_e f=-1\)
 becomes, by~(\ref{polar laplacian}), \(f''(r)+(1/r) f'(r)=-1\) in
 \((a,b)\). The general solution is \(f(r)= -{r}^2/4 + A\log {r} +
 B\), and the boundary conditions \(f(a)=f(b)=0\) make \(A=
 \dfrac{b^2-a^2}{4(\log b-\log a)}\) and \(B= \dfrac{a^2\log b-b^2\log
   a}{4(\log b-\log a)}\).

\subsection{Angular domain}\label{subsec:cone}

Let $V\subset \mathbf{R}^2$ be an angular domain of angle $\alpha$.
In polar coordinates, $V$ is, up to isometry, the set of 
$(r,\theta)$ with $r>0$ and $\theta\in (-\alpha/2,\alpha/2)$.

\begin{svgraybox}
  The expected first exit time from an angular domain \(V\) of angle
  \(\alpha\) is infinite if \(\alpha\ge \pi/2\), and is finite if
  \(\alpha< \pi/2\) and given by
\[ E_{(r,\theta)}[T_V]= \frac{r^2}{4}\left(\frac{\cos 2\theta}{\cos
    \alpha}-1\right).\]
\end{svgraybox}

Because an angular domain is not a relatively compact domain, there is
no guarantee of existence or of uniqueness of solutions to
Equation~\ref{equation1}. The expected first exit time,
\(f=E_\bullet[T_V]\), from \(V\) has two other properties that, if
finite, will characterize it uniquely among the solutions to that
equation.
\begin{enumerate}
\item[(a)]\label{item:a} The expected first exit time function is
  homogeneous of order \(2\).  Indeed, for \(\lambda>0\), the
  mapping \((x,y)\mapsto(\lambda x,\lambda y)\) is a dilation of the
  Euclidean metric that leaves \(V\) invariant. The Euclidean heat
  kernel density at \(\lambda p=(\lambda x, \lambda y)\) at time \(t\)
  deposited in \(\lambda p'=(\lambda x',\lambda y')\) is the
  Euclidean heat kernel density at \(p=(x,y)\) at time \(t/\lambda^2\)
  deposited in \(p'=(x',y')\):
  \[ \frac{1}{2\pi t} e^{-|\lambda p-\lambda p'|^2/4t} d(\lambda
  x')d(\lambda y') = \frac{1}{2\pi (t/\lambda^2)} e^{-|p-
    p'|^2/4(t/\lambda^2)} dx' dy'.\] 
  The effect of such dilation is to
  rescale Brownian motion times by a factor of \(\lambda^2\), and so
  the expected first exit time \(f\) must satisfy
  \[f(\lambda x,\lambda y) = \lambda^2 f(x,y);\]  in polar coordinates
  \[ f(\lambda r,\theta) =\lambda^2 f(r,\theta).\]
\item[(b)]\label{item:b} In polar coordinates, the expected first exit time
  satisfies \(f(r,\theta)=f(r,-\theta)\) because reflection about the
  axis of \(V\) is an isometry that leaves \(V\) invariant.
\end{enumerate}

Consequently, by~(a), the expected first exit time \(f\) is completely
determined by the values \(f(1,\theta)\), and thus it can be written
as \(f(r,\theta)=r^2 h(\theta)\), where \(h\) is positive and
symmetric on $(-\alpha/2,\alpha/2)$. Writing out the differential
equation \(\triangle_e f=-1\) for \(f(r,\theta)=r^2 h(\theta)\) in
polar coordinates~(\ref{polar laplacian}) results in  the following differential equation for
\(h\) on \((-\alpha/2,\alpha/2)\):
\begin{equation}\label{equation:cone}
   h''(\theta) + 4 h(\theta) = -1. 
\end{equation} 
The boundary condition \(f\equiv 0\) on \(\partial C\) results in the boundary
condition \(h(\pm \alpha/2) =0\). Note also that \(h>0\) and that, by~(b),
\(h\) is symmetric about \(0\in [-\alpha/2,\alpha/2]\).
 
The general solution to Equation~(\ref{equation:cone}) is of the form
\[ h(\theta) = A \cos 2\theta + B \sin 2\theta - 1/4. \] The symmetry
of $h$ about $0$ in $[-\alpha/2, \alpha/2]$ implies that $B=0$, and the
initial conditions impose that $ A = 1/(4 \cos \alpha)$. Therefore,
\[ h(\theta) = \frac{\cos 2\theta}{4\cos \alpha}-\frac{1}{4}.\] 
Because this holds for all $\theta$ between $-\alpha/2$ and
$\alpha/2$, and $h>0$, you must have $\cos \alpha >0$ and with the
same sign as $\cos 2\theta$, and so $\alpha<\pi/2$. Writing
\(f(r,\theta)=r^2h(\theta)\) confirms the statement at the beginning
of this section.

In rectangular coordinates, an angular domain of angle \(\alpha<\pi/2\) is
isometric to the domain, \(V\), consisting of all \((x,y)\)
such that \(x>m|y|\), where \(0<m=\tan\alpha/2<1\).

\begin{svgraybox}
The expected first exit time from \(V=\{(x,y)\mid x>m|y|\}\), with \(0<m<1\), is
\begin{equation}\label{angular:exit}
E_{(x,y)}[T_V]= \frac{1}{2-2m^2} \left(m^2x^2-y^2\right).
\end{equation}
Any solution, \(g\), to Equation~(\ref{equation1}) on \(V\) is given by 
\begin{equation*}
g(x,y)= \frac{1}{2-2m^2} \left(m^2x^2-y^2\right) +
C (x^2+y^2)^{\pi/\alpha} \cos \dfrac{\pi\arctan (y/x)}{\alpha},
\end{equation*}
for some constant \(C\ge0\).
\end{svgraybox}

The  second summand in the expression for \(g\)
is justified in a manner similar to that of the case of the parabola at the end 
 Section~\ref{section:parabola}. In this case, you consider the
conformal representation, \(\phi\), of the angular domain \(V\) above onto the
right half plane given by \(\phi(z)=z^{\pi/\alpha}\), and then you
use~(\ref{PoissonKernelRep}) to show that any non-negative harmonic
function on \(V\) that is identically \(0\) on \(\partial V\)
is of the form \(z\mapsto C \Re(z^{\pi/\alpha})\), which is as stated
above.

\subsection{Domain bounded by a hyperbola}

\subsubsection{Convex domain}\label{convex hyperbola}
Let \(D\) be a convex domain in \(\mathbf{R}^2\) bounded by a
hyperbola. Up to isometry, this hyperbola is given by an equation of
the form \((x/a)^2-(y/b)^2=1\), with \( a, b>0\), and \(D\) consists
of all \((x,y)\) such that \((x/a)^2-(y/b)^2>1\) and \(x>0\).

\begin{svgraybox}
If \(b\ge a\), then the expected first exit time from \(D\) is infinite.
\end{svgraybox}

Indeed, \(D\) contains the angular domain \(V\) bounded by the lines \(y=\pm
m (x-a)\) with \(m=b/a\) and \(x>a\), and so the expected first exit times
\(E_{\bullet}[T_D]\ge E_\bullet[T_V]\). If \(b\ge a\), then \(V\) has
angle \(2\arctan m \ge \pi/2\), and so, as established in
Section~\ref{subsec:cone},  the expected first exit time from
\(V\) is infinite.


\begin{svgraybox}
  If \(b<a\), then the expected first exit time from \(D\) is given by the
  function \[g(x,y)=\dfrac{1}{2-2m^2}\left(m^2x^2-y^2-b^2\right),\]
  for all \((x,y)\in D\).\end{svgraybox}

It is plain that \(g\) is a solution to Equation~(\ref{equation1}). If
\(g\) is not the expected first exit time from \(D\), then it is not the
minimal solution to~(\ref{equation1}), and so the expected first exit time
from \(D\) is of the form \(E_\bullet[T_D]=g-u\), where \(u\) is a
positive, harmonic function on \(D\) satisfying \(u\equiv 0\) on
\(\partial D\).

If \(m=b/a<1\), the expected first exit time \(E_\bullet[T_V]\) is finite
and given by~(\ref{angular:exit}) (after a horizontal shift), and
is a minorant for \(E_\bullet[T_D]=g-u\). Thus
\[\dfrac{1}{2-2m^2}\left(m^2(x-a)^2-y^2)\right)\le \dfrac{1}{2-2m^2}\left(m^2x^2-y^2-b^2\right) -u(x,y),\]
or
\begin{equation}\label{hyperbola estimate} 
  u(x,y) <  \frac{b^2}{1-m^2}\left(\frac{x}{a}-1\right),
\end{equation}
for all \((x,y)\) in \(V\).

You will now reach a contradiction as follows. The function \(\phi\)
given by (appropriate branches taken) 
\begin{align*}
\phi(z) & =c \cosh\left(\dfrac{2\mu}{\pi}
  \operatorname{arcosh} z\right)\\ & =\frac{c}{2}\left((z+\sqrt{z^2-1})^{2\mu/\pi}+(z+\sqrt{z^2-1})^{-2\mu/\pi}\right),
\end{align*}
 where \(\mu=\arctan m\) and \(c=\sqrt{a^2+b^2}\) the linear eccentricity, is a
conformal representation of the right half plane \(\{x=\Re z>0\}\)
onto \(D\) which takes the boundary \(\{x=\Re z= 0\}\) onto \(\partial
D\) and the ray \(\{x \ge 0\}\) onto the ray \( \{x\ge a\}\). (As an
aid in visualizing this mapping, you recall that \(\cosh\) takes the
horizontal line through \(\beta i\) \((0<\beta<\pi/2)\) onto the
right branch of the hyperbola of equation \(x^2/\cos^2\beta -
y^2/\sin^2\beta=1\), cf.~\cite[3.4.2]{Ahlfors} or \cite{Kober} for more background.)

Then, on the one hand, you deduce from inequality~(\ref{hyperbola
  estimate}) that
\[ u(\phi(x)) < \frac{b^2}{1-m^2}\left(\frac{\phi(x)}{a}-1\right),\]
because \(\phi(z)\) is real for \(z\) real, and from this inequality
that \(\displaystyle\lim_{x\to \infty} \dfrac{u(\phi(x))}{x}=0\),
because, for real \(x>1\),
\(\phi(x)=\dfrac{c}{2}\left((x+\sqrt{x^2-1})^{2\mu/\pi}+(x+\sqrt{x^2-1})^{-2 \mu/\pi}\right)\), hence
  \(\phi(x)=O(x^{2\mu/\pi})\),  and so \(\displaystyle
  \lim_{x\to \infty}\dfrac{\phi(x)}{x}=0\) because \(\mu=\arctan m <\pi/2\).

  On the other hand, the composite function \(u\circ \phi\) is
  positive and harmonic on \(\Re z>0\) and is identically \(0\) on the
  boundary \(\Re z=0\). Therefore, by~(\ref{PoissonKernelRep}),
  \(u\circ \phi(x,y)=C x\), for some constant \(C>0\), and so
  \(\displaystyle\lim_{x\to \infty} \dfrac{u(\phi(x))}{x}=C>0\).


\subsubsection{Concave domain}

A (concave) domain bounded by the two branches of a hyperbola is
isometric to the domain \(D=\{ (x,y)\in \mathbf{R}^2\mid
x^2/a^2-y^2/b^2>-1\}\), for some \(a, b>0\). Let \(m=b/a\) and
\(\mu=\arctan m\). If \(m\ge1\), then the expected first exit time
from \(D\) is infinite because \(D\) contains an angular domain
of angle \(2\mu\ge\pi/2\).

\begin{svgraybox}
  If \(b<a\), the expected first exit time from \(D\) is given by the
  function 
\begin{equation}
\label{exit concave hyperbola} g(x,y)=\dfrac{1}{2-2m^2}\left(m^2x^2-y^2+b^2\right),\end{equation}
  for all \((x,y)\in D\). 
\end{svgraybox}

Certainly, \(g\) is a solution to Equation~(\ref{equation1}) on \(D\),
so \(g=E_{(x,y)}[T_D]+v\), where \(v\) is \(\ge 0\) and harmonic
on \(D\) and extends continuously to \(0\) on \(\partial D\).

The domain \(D\) contains the angular domains \(V_{-}=\{x<0
\text{~\&~} b^2x^2>a^2y^2\}\) and \(V_+ =\{ x>0 \text{~\&~}
b^2x^2>a^2y^2\}\). These domains have angle \(2\mu<\pi/2\), so their
expected first exit time is finite and given by
\(f_{\pm}(x,y)=\dfrac{a^2b^2}{2a^2-2b^2}\left(b^2x^2-a^2y^2\right)\) (same
expression, different domain). By comparison,
\[ \dfrac{a^2b^2}{2a^2-2b^2}\left(b^2x^2-a^2y^2\right) \le
\dfrac{a^2b^2}{2a^2-2b^2}\left(b^2x^2-a^2y^2+1\right)-v(x,y),\]
or,
\begin{equation}\label{v is bounded}
v(x,y)\le \frac{a^2b^2}{2a^2-2b^2},\end{equation} 
for all \((x,y)\in V_{-} \cup
V_{+}\). In particular, \(v\) is bounded on the real axis \(\{y=0\}\).

The mapping \(\psi\) given by 
\begin{align*}
\psi(z) & =\frac{c}{2}\left(z^{2\mu/\pi}
  - \frac{1}{z^{2\mu/\pi}}\right)\\
& =c \sinh \left(\log z^{2\mu/\pi}\right),
\end{align*}
 where \(c=\sqrt{a^2+b^2}\), is a conformal representation of
the right half plane \(\Re z>0\) onto \(D\) that takes the positive
imaginary axis \(\Im z>0\) to the upper branch of the hyperbola
bounding \(D\), takes the negative imaginary axis to the lower branch
of that hyperbola, and takes the positive real axis onto the real
axis (cf.~\cite[3.4.2]{Ahlfors} or~\cite{Kober} for help on constructing this mapping). Therefore, the composite function \(v\circ\psi\) is a positive,
harmonic function on \(\Re z>0\) and that extends continuously by
\(0\) to \(\Im z\ne 0\) on the right half plane. Therefore,
by~(\ref{PoissonKernelRep}), \(v\circ \psi\) is of the form 
\[v\circ \psi(z) = A\Re z + B \Re (1/z)\] for some constants \(A,B\ge
0\). Because of (\ref{v is bounded}), the composite \(v\circ \psi\)
must be bounded on the positive real axis \(\Im z=0\), and therefore
both \(A=0\) and \(B=0\), that is, \(v\equiv 0\), which shows that the
expected first exit time \(E_{\bullet}[T_D]=g\) as was stated at~(\ref{exit concave hyperbola}).

\section{Domains in the Hyperbolic Plane}

The unit disk model for the hyperbolic plane is realized by the unit
disk \(\{x^2+y^2<1\}\) in \(\mathbf{R}^2\), endowed with the metric
\(4{(1-x^2-y^2)^{-2}}\left( dx\otimes dx + dy \otimes dy\right)
\).

In  geodesic polar coordinates \((r,\theta)\) about a point, the hyperbolic metric is 
\( dr\otimes dr + \sinh^2 r d\theta \otimes d\theta\), 
and the  corresponding Laplacian is
\begin{equation}
\label{HyperbolicLaplacia:gpolar}
\triangle_h f = f_{rr} +2\operatorname{cotanh r} f_r +
f_{\theta\theta}.
\end{equation}  

The right half plane model for the hyperbolic plane is realized on
the domain \(\{ x>0\}\)  with the metric
\({x^{-2}}(dx\otimes dx + dy\otimes dy)\), with  Laplacian
 \begin{equation}\label{HyperbolicLaplacian:RightPlane}
\triangle_h f= x^2(f_{xx}+f_{yy}).\end{equation} 
In {E}uclidean polar coordinates \((r, \theta)\), with \(r>0\) and
\(-\pi/2<\theta<\pi/2\), the Laplacian is 
\begin{equation}\label{HyperbolicLaplacian:epolar}
 \triangle_h f= (\cos^2 \theta) \left(r^2 f_{rr} + r f_r + f_{\theta\theta}\right).\end{equation}


\subsection{Hyperbolic disks}


Let \(D\) be a disk of radius \(R\) in the hyperbolic plane. In
geodesic coordinates \((r,\theta)\) based at the center of \(D\), the
Laplacian is \(\triangle_h f=f_{rr}+2\operatorname{cotanh r}f_r +
f_{\theta \theta}.\) The expected first exit time from \(D\) is invariant by
(hyperbolic) rotations about the center of \(D\), therefore \(f=f(r)\)
is a function of the distance \(r\) only (the hyperbolic distance to
the center of the disk), and so the equation \(\triangle_h f=-1\) on
\(D\) becomes, by~(\ref{HyperbolicLaplacia:gpolar}),
\[f''(r)+ 2 \operatorname{cotanh}(r) f'(r) =-1\] on \((0,R)\), with
the boundary conditions that \(\displaystyle \lim_{r\to 0} f(r)\)
exists and that \(f(R)=0\). The general solution is
\(f(r)= -(r/2)\operatorname{cotanh} r + A\operatorname{cotanh}r+B
\), and the boundary conditions imply that \(A=0\) and that 
\(B=(R/2) \operatorname{cotanh}R\). 

\begin{svgraybox}
The expected first exit time from a hyperbolic disk of radius \(R\) is
\begin{equation}\label{exit:hyperbolic disk}
f(r,\theta)=  -\frac{r}{2}\operatorname{cotanh} r +
\frac{R}{2} \operatorname{cotanh} R.\end{equation}
\end{svgraybox}



\subsection{Horodisks}

A horodisk is a domain in the hyperbolic plane that may be
{visualized} as a hyperbolic disk with center at a point on the ideal
boundary of the hyperbolic plane. Any horodisk is isometric to the
domain, \(D\), in the right half plane model consisting of the points
\((x,y)\) with \(x>R\), for some constant \(R>0\).  
\begin{svgraybox}
The expected
first exit time from the horodisk \(D=\{ x>R\}\) is
given by
\begin{equation}
\label{exit:horodisk}
E_{(x,y)}[T_D] = \log\frac{x}{R}.\end{equation}
\end{svgraybox}
The expected first exit time, \(f=E_\bullet[T_D]\), from \(D\) must be
invariant under vertical translations because these are hyperbolic
isometries that leave \(D\) invariant; that is, \(f(x,y+t)=f(x,y)\)
for all \(t\), or \(f(x,y)=f(x)\) for all \(y\). Given this, the
equation \(\triangle_h f=-1\), with \(f\equiv 0\) on
\(x=R\), reduces to
\[f''(x) = -x^{-2}\] on \((R,\infty)\), with $f(R) = 0$, that is
\(f(x,y)= \log\dfrac{x}{R}\), as advertised at~(\ref{exit:horodisk}).

\subsection{Neighborhoods of geodesics}\label{section:geodesicNhood}

A neighborhood, \(D\), of a geodesic in the hyperbolic plane is
determined, up to isometry, by its radius: \(D=D(R)\). In the right
half plane model and in {E}uclidean polar coordinates \((r,\theta)\) as
in~(\ref{HyperbolicLaplacian:epolar}), the geodesic is \(\theta=0\), and
\(D(R)\) is the set of all \((r,\theta)\) with \(r>0\) and
\(-\alpha<\theta<\alpha\), where \(\log\dfrac{\cos\alpha}{1-\sin\alpha}=R\).


Because the mapping \((r,\theta)\mapsto (\lambda r, \theta)\) is an
isometry of the hyperbolic plane that preserves \(D\) (it is a
hyperbolic translation along the geodesic \(\theta=0\)), the expected
first exit time, \(f=E_\bullet[T_D]\), must satisfy \(f(r,\theta)=f(\lambda
r, \theta)\) for all \(\lambda>0\). That is, \(f(r,\theta)=f(\theta)\)
is a function of \(\theta\) only.  So writing out the equation \(\triangle_h f=-1\)  in coordinates
\((r,\theta)\), results in, cf.~(\ref{HyperbolicLaplacian:epolar}),
\[ f''(\theta) = \dfrac{-1}{\cos^2\theta}.\] The general solution is
\(f(\theta)=\log\cos \theta + A\theta+B\).  The expected first exit time is
invariant under reflection about \(\theta=0\), so the constant
\(A=0\), and the boundary condition \(f(\alpha)=0\) forces \(B=-\log
\cos\alpha\).

\begin{svgraybox} The expected first exit time from \(D=D(\alpha)\) is given
  by
\[
f{(r,\theta)}= \log\dfrac{\cos \theta}{\cos\alpha}.
\]
Any solution to Equation~(\ref{equation1}) on \(D\) is of the form
\[f(r,\theta)+ \left(A r^{\pi/2\theta}+ \frac{B}{r
^{\pi/2\theta}}\right)\cos \frac{\pi \theta}{2\alpha},
\]
for some constants \(A,B\ge0\).
\end{svgraybox}

Indeed, the composite of
a non-negative harmonic function, \(u\), on \(D\) that is 
\(0\) on the boundary  \(\partial D\) and the conformal representation of
the right half plane \(\{ x>0\}\) onto \(D\) given by
\(\phi:(r,\theta)\mapsto (r^{2\alpha/\pi}, 2\alpha/\pi\theta)\) is a
non-negative harmonic function, \(u\circ \phi\), on the right half
plane that is identically \(0\) on the boundary \(\theta=\pm \pi/2\),
except perhaps at the origin. By~(\ref{PoissonKernelRep}), such
function is of the form \((r,\theta)\mapsto (Ar+B/r)\cos\theta\), for some constants
\(A,B\ge 0\).


\subsection{Neighborhood of ideal point}

You may adjust the calculation in the previous section to cover the
case of a one-sided neighborhood of a geodesic. In Euclidean polar
{coordinates} as above, the one-sided neighborhood of radius \(R\) is
isometric to the domain \(D=D(R)\) of all \((r,\theta)\) with \(r>0\)
and \(0<\theta<\alpha(R)\), with \(\alpha(R)\) as in the first
paragraph of Section~\ref{section:geodesicNhood}. The expected first
exit time from \(D(R)\) is
\[E_{(r,\theta)}[T_{D(R)}]= \log\cos \theta -\dfrac{\theta}{\alpha}\log\cos
\alpha.\]

A neighborhood of an ideal point is isometric to the domain \(D(R)\)
when \(R=\infty\), or \(\alpha=\pi/2\). Hence, the expected
first exit time from such domain is infinite because there is no
positive solution to
the differential equation \(f''(\theta)=-1/\cos^2\theta\) on
\((0,\pi/2)\) with boundary conditions \(f(0)=f(\pi/2)=0\).

\bibliographystyle{spbasic}

\end{document}